\documentclass[a4paper, twocolumn]{article}

\usepackage[utf8]{inputenc}
\usepackage[T1]{fontenc}
\usepackage{textcomp}
\usepackage{amsmath, amssymb, amsthm}
\usepackage[margin=2cm]{geometry}
\usepackage{tikz-cd}
\usepackage{mdframed}
\usepackage{microtype}
\usepackage{graphicx,tipa}
\newcommand{\arc}[1]{{
	\setbox9=\hbox{\ensuremath{#1}}
  \ooalign{\resizebox{\wd9}{\height}{\texttoptiebar{\phantom{A}}}\cr{\ensuremath {#1}}}}}
\usepackage{enumitem}

\usepackage{import}
\usepackage{xifthen}
\usepackage{todonotes}
\usepackage{pdfpages}
\usepackage{transparent}

\usepackage[backend=biber]{biblatex}
\usepackage{hyperref}
\usepackage{cleveref}

\newcommand{\sphere}{\mathbb{S}^2}
\newcommand{\E}{\mathbb{E}}

\newcommand{\mink}{\mathbb{M}^3}

\DeclareMathOperator{\area}{Area}
\DeclareMathOperator{\der}{d}

\theoremstyle{definition}

\newtheorem{theorem}{Theorem}
\newtheorem*{remark}{Remark}

\author{Micha\"el Maex}
\date{}
\title{A synthetic proof of the spherical and hyperbolic Pythagorean theorem on models in Euclidean and Minkowski space}

\begin{document}

\maketitle

This paper is a response to Paolo Maraner's previous paper \cite{maranerSphericalPythagoreanTheorem2010} in \emph{The Mathematical Intelligencer} and the subsequent letter to the editors by Victor Pambuccian \cite{pambuccianMariaTeresaCalapsos2010}.
These papers discuss a version of the Pythagorean theorem that not only holds in flat Euclidean space, but also in spherical and hyperbolic geometry. The hyperbolic case was first discovered by Maria Teresa Calapso \cite{calapso}, and rediscovered independently and generalized to the spherical case by Paolo Maraner \cite{maranerSphericalPythagoreanTheorem2010}. 

It seems that all known proofs are analytic in nature, giving very little geometric intuition. 
A synthetic proof is still missing, which is the final question of \cite{maranerSphericalPythagoreanTheorem2010}.
It was frustrating to me that such a simple generalisation of such a famous theorem does not have a proof that is more geometric in nature. 
In this paper, we will present a proof that is much more geometric in nature than previously known proofs.

The paper has two parts. The first part, \cref{sec:the_theorem,sec:a_simple_proof_on_the_sphere}, is aimed at any enthusiastic reader with a high school level understanding of geometry.
Here, we will introduce and motivate the statement of the theorem and go over its history. We then give a proof without worrying about precise axiomatic underpinnings.
An essential ingredient is considering spherical space as the surface of a sphere in 3-dimensional Euclidean space $\mathbb E^{3}$. 

The second part, \cref{sec:the_hyperbolic_case,sec:axioms_and_embeddings}, is for geometry enthusiasts, either professional mathematicians or casual enjoyers of unusual geometries. 
In this part we will discuss how the proof actually generalises to hyperbolic space as well, once we change from considering the sphere in 3-dimensional Euclidean space to the hyperboloid in Minkowski 3-space. 
At the end we muse on which axioms/postulates are needed to phrase the theorem and carry out the proof.

\section{The theorem and its history} \label{sec:the_theorem}

We are all familiar with the classical Pythagorean theorem in the plane. 
\begin{theorem}[Original Pythagorean theorem]
	Let $ABC$ be a right angled triangle with right angle $A$ in the Euclidean plane, $\mathbb E^2$.
	Then the squares on the sides $AB$ and $AC$ together have the same area as the square on the side $BC$.
\end{theorem}

\begin{figure}
	\centering
	\def\svgwidth{\columnwidth}
	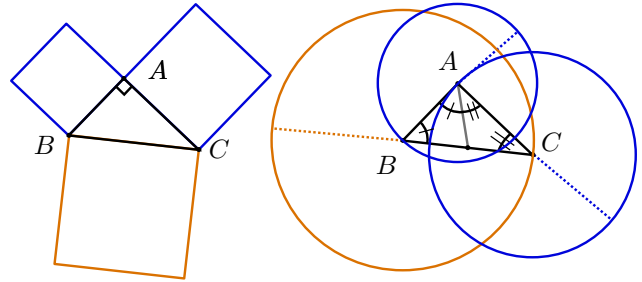

	\caption{The original (left) and reformulated (right) Pythagorean theorem in $\mathbb E^2$. In either formulation the area of the blue figures sum to the area of the yellow figure.}
	\label{fig:comparison_original}
\end{figure}

In fact, we can replace the squares on the sides with any shape with whose dimensions are proportional to the side, such as disks with the sides as radii. 
A right angled triangle can also be described as a triangle where one angle is equal to the sum of the two other angles. We call such a triangle \emph{properly angled}, and its largest angle the \emph{proper angle}. 
So we may also phrase the Pythagorean theorem as:
\begin{theorem}[Absolute Pythagorean theorem]
	Let $ABC$ be a properly angled triangle with proper angle $A$. 
	Write $O_{AB}, O_{BC}, O_{CA}$ for the disks with radii $|AB|, |BC|, |CA|$ respectively. 
	Then 
	\[
		\area (O_{BC}) = \area (O_{AB}) + \area( O_{CA})
	.\] 
\end{theorem}
See \cref{fig:comparison_original} for a comparison between the two formulations. 

Notice that I didn't write that the triangle is taken in the Euclidean plane $\E^2$. 
This is because this version of the theorem is actually true in spherical and hyperbolic geometry too!
We say that a theorem that holds in Euclidean, spherical and hyperbolic geometry is true in \emph{absolute} geometry. 
The hyperbolic case was first shown by Maria Teresa Calapso in \cite{calapso}, and was discovered independently and extended to the spherical case by Paolo Maraner in \cite{maranerSphericalPythagoreanTheorem2010}. 
To the best of my knowledge all known proofs are analytic in nature, involving formulae for the sides of triangles and areas of disks.

Let $ABDC$ be a quadrilateral. 
We say that $ABDC$ is \emph{equiangular} if all four angles are equal. 
These are precisely quadrilaterals for which $\arc {AD}, \arc{BC}$ are equally long and intersect each other in their midpoints, 
in particular the points $A,B, C, D$ lie on a circle with center at the intersection of the diagonals. 
It is known that every proper triangle is half of an equiangular quadrilateral cut along the diagonal e.g.\ see \cite{maranerSphericalPythagoreanTheorem2010}. 

\begin{figure}[h]
	\centering
	\def\svgwidth{\columnwidth}
\begingroup%
  \makeatletter%
  \providecommand\color[2][]{%
    \errmessage{(Inkscape) Color is used for the text in Inkscape, but the package 'color.sty' is not loaded}%
    \renewcommand\color[2][]{}%
  }%
  \providecommand\transparent[1]{%
    \errmessage{(Inkscape) Transparency is used (non-zero) for the text in Inkscape, but the package 'transparent.sty' is not loaded}%
    \renewcommand\transparent[1]{}%
  }%
  \providecommand\rotatebox[2]{#2}%
  \newcommand*\fsize{\dimexpr\f@size pt\relax}%
  \newcommand*\lineheight[1]{\fontsize{\fsize}{#1\fsize}\selectfont}%
  \ifx\svgwidth\undefined%
    \setlength{\unitlength}{693.15092793bp}%
    \ifx\svgscale\undefined%
      \relax%
    \else%
      \setlength{\unitlength}{\unitlength * \real{\svgscale}}%
    \fi%
  \else%
    \setlength{\unitlength}{\svgwidth}%
  \fi%
  \global\let\svgwidth\undefined%
  \global\let\svgscale\undefined%
  \makeatother%
  \begin{picture}(1,0.44949016)%
    \lineheight{1}%
    \setlength\tabcolsep{0pt}%
    \put(0,0){\includegraphics[width=\unitlength,page=1]{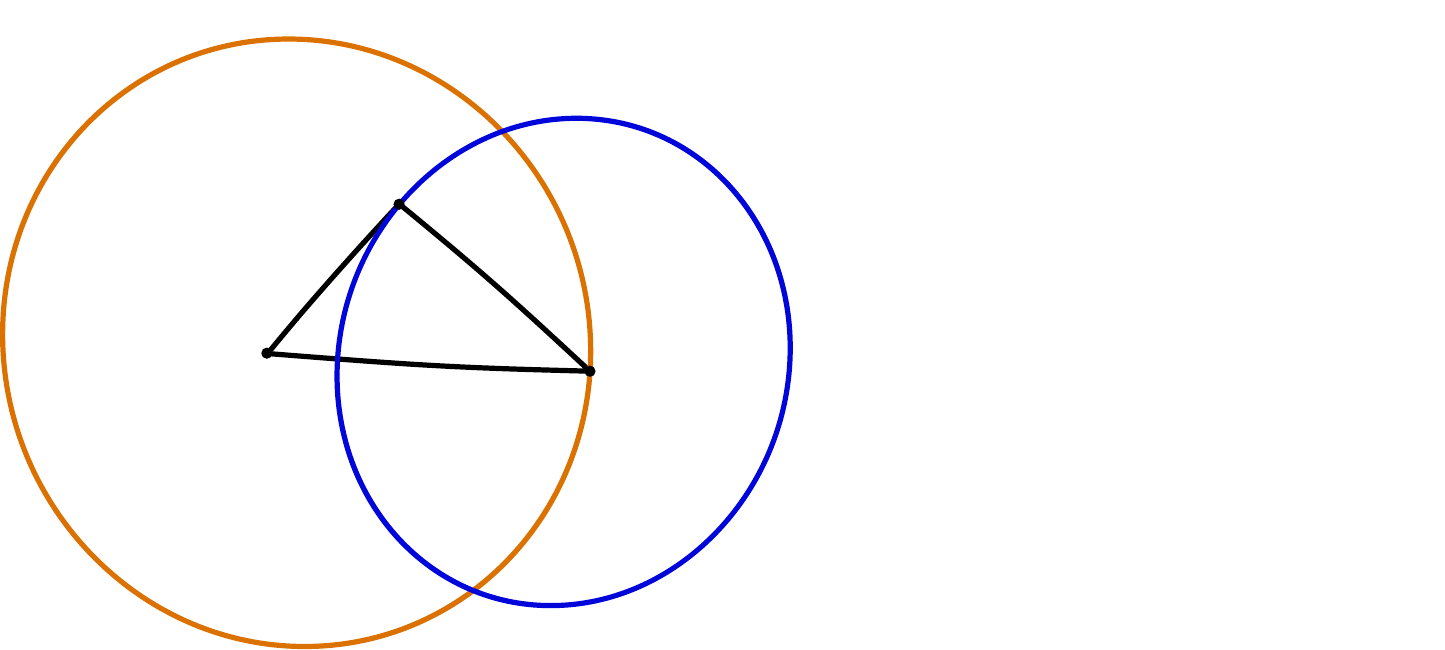}}%
    \put(0.25073883,0.32141809){\makebox(0,0)[lt]{\lineheight{0}\smash{\begin{tabular}[t]{l}$A$\end{tabular}}}}%
    \put(0.17040856,0.22823366){\makebox(0,0)[lt]{\lineheight{0}\smash{\begin{tabular}[t]{l}$B$\end{tabular}}}}%
    \put(0.42365264,0.17621148){\makebox(0,0)[lt]{\lineheight{0}\smash{\begin{tabular}[t]{l}$C$\end{tabular}}}}%
    \put(0,0){\includegraphics[width=\unitlength,page=2]{pythagorean.pdf}}%
    \put(0.83510793,0.01741151){\makebox(0,0)[lt]{\lineheight{0}\smash{\begin{tabular}[t]{l}$D$\end{tabular}}}}%
    \put(0.75699445,0.26341905){\makebox(0,0)[lt]{\lineheight{0}\smash{\begin{tabular}[t]{l}$A$\end{tabular}}}}%
    \put(0.67977677,0.18103354){\makebox(0,0)[lt]{\lineheight{0}\smash{\begin{tabular}[t]{l}$B$\end{tabular}}}}%
    \put(0.9335239,0.13267485){\makebox(0,0)[lt]{\lineheight{0}\smash{\begin{tabular}[t]{l}$C$\end{tabular}}}}%
    \put(0,0){\includegraphics[width=\unitlength,page=3]{pythagorean.pdf}}%
  \end{picture}%
\endgroup%

	\caption{The diagonal reformulation of the absolute Pythagorean theorem}
	\label{fig:diagonal_pythagorean}
\end{figure}

The version of the Pythagorean theorem that we will prove is a reformulation in terms of the diagonal of an equiangular quadrilateral. In this case we can take all the involved circles to be concentric, which works out nicely in the proof. 
\begin{theorem}[Diagonal Pythagorean theorem]\label{thm:diag_pyth}	
	Let $ABDC$ be an equiangular quadrilateral in Euclidean, spherical or hyperbolic space. 
	Let $O_{AB}, O_{AC}, O_{AD}$ be the disks with center $A$ and radii $|AB|, |AC|, |AD|$ respectively. 
	Then \[
		\area(O_{AB}) + \area(O_{AC}) = \area(O_{AD})
	.\] 
\end{theorem}
I'll leave it as an exercise to the motivated reader to see why this reformulation is equivalent to the previous one.
\Cref{fig:diagonal_pythagorean} could be a hint.

\section{A geometric proof on the sphere} \label{sec:a_simple_proof_on_the_sphere}

\subsection{What is a sphere?} \label{sec:what_is_spherical_geometry}

This may seem like a silly question. 
After all, we all have some intuitive understanding of spheres from playing football or any other ball sport.
If I ask you to picture a sphere, you will probably imagine the surface of some 3-dimensional object. 
It might be surprising to learn that mathematicians can think about spheres without needing to imagine them inside any surrounding space. 
In axiomatic geometry for example, we could take the normal axioms of flat Euclidean space but replace the parallel postulate by another axiom.
If this sounds complicated, don't worry. 
In this paper, we will only consider the sphere as existing in 3-dimensional space. In fact, this is vital to our argument!

Euclidean space $\mathbb E^3$ consists of triplets $(x, y, z)$.
At school we learn that a sphere consists of all points lying at a fixed distance to a center point. 
If we choose the center at $(0,0,0)$ and the radius $r > 0$ then our model of the sphere is the surface given by all triplets satisfying the equation \[
x^2 + y^2 + z^2 = r^2
.\] 
These are simply all points at distance $r$ from the center point $(0,0,0)$.
Let's denote our sphere by $\sphere$ and the surrounding 3-dimensional Euclidean space by $\E^3$. 

The ``straight'' lines without ends on the sphere are the intersections of $\sphere$  with planes in $\E^3$ that pass through the center.
We call these lines on sphere \emph{great circles}, e.g.\ the equator on Earth. 
Parts of great circles that are bounded between two points are called \emph{arcs}.
These are analogous to line segments in Euclidean geometry.
If $A, B$ are two points $\sphere$ that are not antipodal they will lie on exactly one great circle. 
We will write  $\arc{AB}$ for the shortest arc on this great circle bounded by $A,B$, and we write $\overline{AB}$ for the corresponding straight line segment in $\E^3$.

Just like in Euclidean geometry, circles on $\sphere$ are all points that lie at a fixed distance from some center point. 
It turns out the circles on $\sphere$ are precisely the non-empty intersections of $\sphere$ and some plane (not necessarily through the center) in $\E^3$.
For example, the lines of latitude on earth are circles with the north/south pole as center, and they also cut out by planes parallel to the equator.

\subsection{The bread-crust theorem} \label{sec:the_bread_crust_theorem}
We need a way to measure the area of a disk on the sphere. 
Fortunately, there is an easy way to do this if we think of the sphere as a surface in 3-dimensional space.
Let $\sigma$ be the slice of $\sphere$ cut out by two parallel planes, each intersecting $\sphere$ and having a distance $d$ apart. 
Then 
\[
	\area(\sigma) = d \cdot 2\pi r
.\] 
This relation is originally due to Archimedes, and those who are interested in learning more about this can look at the beautiful explanation by the wonderful science communicator Grant Sanderson on his YouTube channel 3Blue1Brown \cite{3blue1brownWhySpheresSurface2018}.

I like to refer to this theorem as the bread-crust theorem, because it means that if you take a slice of a perfectly spherical loaf of bread, the amount of crust on the slice is proportional to the width of the slice. 
In particular, two slices of the same width will always have the same amount of crust, regardless of whether they are taken from the middle or the edge of the loaf.

\subsection{The actual proof} \label{sec:the_actual_proof}

Start with an equiangular quadrilateral $ABDC$ in $\sphere$.
For the sake of simplicity, we may assume that  $A$ lies at the north pole $(0,0,r)$. 
As we have discussed in \cref{sec:the_theorem}, the vertices of an equiangular quadrilateral $ABDC$ lie on some circle, which we will call $O_{ABDC}$.

\begin{figure}[h]
	\centering
	\def\svgwidth{\columnwidth}
\begingroup%
  \makeatletter%
  \providecommand\color[2][]{%
    \errmessage{(Inkscape) Color is used for the text in Inkscape, but the package 'color.sty' is not loaded}%
    \renewcommand\color[2][]{}%
  }%
  \providecommand\transparent[1]{%
    \errmessage{(Inkscape) Transparency is used (non-zero) for the text in Inkscape, but the package 'transparent.sty' is not loaded}%
    \renewcommand\transparent[1]{}%
  }%
  \providecommand\rotatebox[2]{#2}%
  \newcommand*\fsize{\dimexpr\f@size pt\relax}%
  \newcommand*\lineheight[1]{\fontsize{\fsize}{#1\fsize}\selectfont}%
  \ifx\svgwidth\undefined%
    \setlength{\unitlength}{585.14866146bp}%
    \ifx\svgscale\undefined%
      \relax%
    \else%
      \setlength{\unitlength}{\unitlength * \real{\svgscale}}%
    \fi%
  \else%
    \setlength{\unitlength}{\svgwidth}%
  \fi%
  \global\let\svgwidth\undefined%
  \global\let\svgscale\undefined%
  \makeatother%
  \begin{picture}(1,0.53188537)%
    \lineheight{1}%
    \setlength\tabcolsep{0pt}%
    \put(0,0){\includegraphics[width=\unitlength,page=1]{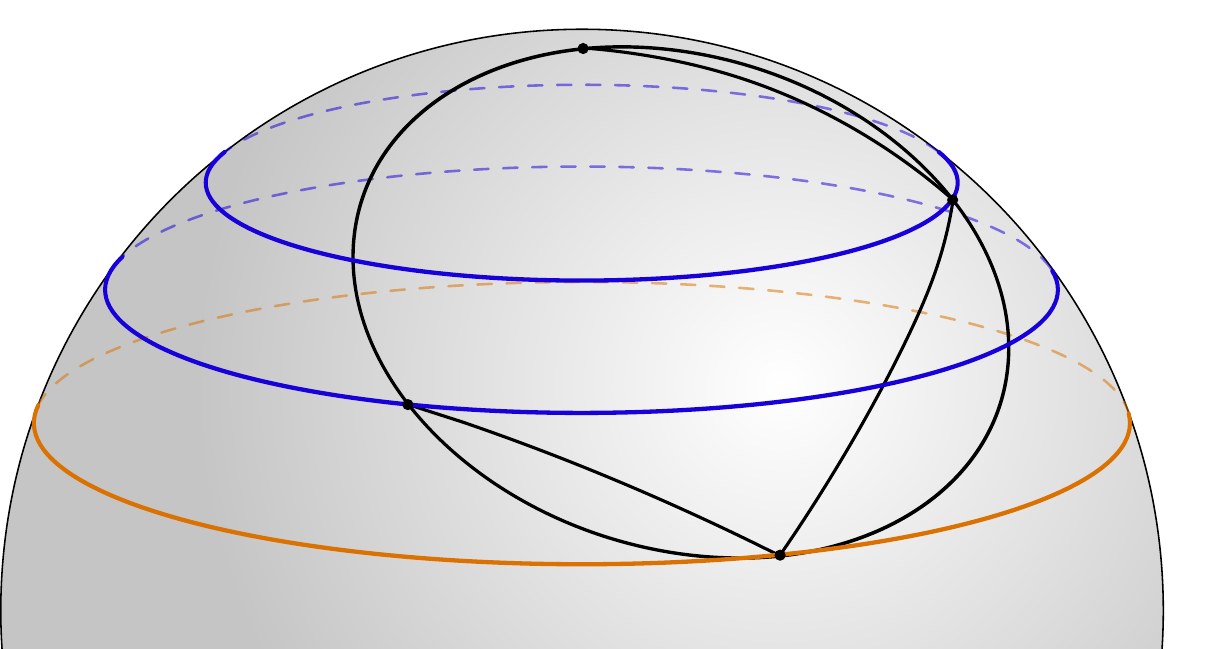}}%
    \put(0.6485407,0.03307539){\makebox(0,0)[lt]{\lineheight{0}\smash{\begin{tabular}[t]{l}$D$\end{tabular}}}}%
    \put(0.47068076,0.52149802){\makebox(0,0)[lt]{\lineheight{0}\smash{\begin{tabular}[t]{l}$A$\end{tabular}}}}%
    \put(0.28582375,0.15152909){\makebox(0,0)[lt]{\lineheight{0}\smash{\begin{tabular}[t]{l}$B$\end{tabular}}}}%
    \put(0.69915951,0.36002883){\makebox(0,0)[lt]{\lineheight{0}\smash{\begin{tabular}[t]{l}$C$\end{tabular}}}}%
    \put(0.77818775,0.41615471){\makebox(0,0)[lt]{\lineheight{0}\smash{\begin{tabular}[t]{l}$O_{AC}$\end{tabular}}}}%
    \put(0.86853912,0.3230096){\makebox(0,0)[lt]{\lineheight{0}\smash{\begin{tabular}[t]{l}$O_{AB}$\end{tabular}}}}%
    \put(0.92886401,0.20501826){\makebox(0,0)[lt]{\lineheight{0}\smash{\begin{tabular}[t]{l}$O_{AD}$\end{tabular}}}}%
    \put(0,0){\includegraphics[width=\unitlength,page=2]{spherical.pdf}}%
  \end{picture}%
\endgroup%

	\caption{The equiangular quadrilateral $ABDC$ and the circles from the spherical Pythagorean theorem}
	\label{fig:sphere}
\end{figure}
Any circle of $\sphere$ is the intersection of the sphere with a plane in $\E^{3}$.
So, let $P$ the plane such that $O_{ABDC} = P \cap \sphere$.
We know that opposite sides of an equiangular quadrilateral are congruent arcs in $\sphere$ and thus they are also congruent in $\E^{3}$. 
In particular, the distances in $\E^3$ between the end points are the same, i.e.\ $|AB| = |CD|$ and  $|AC| = |BD|$.
Thus in $P$ the quadrilateral $ABDC$ is a parallelogram. 
Therefore, we have the following equality of vectors \[
\overrightarrow{BA} + \overrightarrow{CA} = \overrightarrow{DA}
.\]

As the disks $O_{AB}, O_{AC}, O_{AD}$ have the north pole as center, their boundary circles are cut out by planes parallel to the  $xy$-plane. 
As a slice, these disks are cut out by these parallel planes and the parallel plane touching the north pole.
Thus the thickness of the slices is given by $z$-components of the vectors above. 
Using the bread-crust theorem we see that
\begin{align*}
	\area(O_{AD}) &= \overrightarrow{DA}_z \cdot 2\pi r \\
		      &=\left( \overrightarrow{BA}_z + \overrightarrow{CA}_z \right)\cdot 2\pi r  \\
		      &= \area(O_{AB}) + \area(O_{AC}) 
,\end{align*}
which finishes the proof of \Cref{thm:diag_pyth} in the spherical case.
\begin{remark}
	In $\E^3$ the quadrilateral  $ABDC$ is not only a parallelogram but in fact a rectangle.
	So, the proper triangles and equiangular rectangle on the sphere are precisely the points that in $\E^3$ form a right angled triangle or rectangle respectively.
\end{remark}

\section{The hyperbolic case} \label{sec:the_hyperbolic_case}

In order to carry out the proof in \cref{sec:the_actual_proof}, we only need the following properties of the embedding of $\sphere$ in $\E^3$.
\begin{enumerate}[label= (\alph*)]
	\item Isometries of $\sphere$ are the restrictions of isometries of $\E^3$. 
		This is used, for example, when we want to conclude that if $\arc{AB}$ is congruent to $\arc{CD}$ then $|AB| = |CD|$.
	\item The circles in $\sphere$ are precisely the non-empty intersections with a Euclidean plane in $\E^3$.
	\item The bread-crust theorem holds for a slice that is cut out by two parallel Euclidean planes.
\end{enumerate}

The sphere has constant curvature $r^2$. 
So we can vary the radius of $\sphere$ to get a surface with any constant positive curvature. 
Suppose that somehow we could get a sphere with imaginary radius, making $r^2$ negative. 
Then this would give a surface with constant negative curvature, i.e.\ the hyperbolic plane. 

This seemingly nonsensical idea can be made precise by changing the metric $\E^3$ for the indefinite metric \[
x^2 + y^2 - z^2
.\] 
This is known as Minkowski $3$-space, which we will denote by $\mink$.

We can look at all points at imaginary distance $i\cdot r$ from the origin, i.e.\ the solutions to the equation  \[
x^2 + y^2 - z^2 = - r^2
.\] 
This consists of two connected components. Let $\mathbb{H}^2$ be the connected component with $z > 0 $.
This is known as the hyperboloid model of hyperbolic geometry. See \cite[§ 2.3]{thurstonThreedimensionalGeometryTopology1997}\cite{reynoldsHyperbolicGeometryHyperboloid1993} for more details.  

It turns out that the embedding of $\mathbb{H}^2$ into $\mink$ also satisfies the properties mentioned above, which allow the proof to work. 
\begin{enumerate}[label= (\alph*)]
	\item The isometries of $\mathbb{H}^2$ are precisely the restrictions of the linear isometries of $\mink$ that map $\mathbb{H}^2$ to itself (instead of $- \mathbb{H}^2$) \cite[sec 5]{reynoldsHyperbolicGeometryHyperboloid1993}.
	\item Let $O$ be a circle with center $c$. 
		By the previous point, we may assume that  $c$ is the point $(0, 0, r)$.
		Then it is clear that $O$ is the intersection of $\mathbb{H}^2$ with some plane, $P$ parallel to the $xy$-plane. 
		It is also clear that the Minkowski metric restricts on $P$ to the usual Euclidean metric on it. 
	\item Again by using symmetry, we may assume that the slice $\sigma$ is cut by two planes parallel to the $xy$-plane, say $P_1: z = z_1, P_2: z = z_2$ with $r\le z_1 \le z_2$. 
		The hyperboloid is the surface of revolution around the $z$-axis of the curve $x = f(z) = \sqrt{z^2 - r^2}$. 
		The usual area formula for a revolution surface, but derived with the Minkowski metric yields 
		\begin{align*}
		\area(\sigma) &= 2\pi\int_{z_1}^{z_2} f(z)\cdot \sqrt{(f'(z))^2 -1}  \der z \\
			      &= 2\pi \int _{z_1}^{z_2}\sqrt{z^2 - r^2} \sqrt{\frac{z^2}{z^2 - r^2} - 1}   \der z\\
			      &=  2\pi \int_{z_1}^{z_2} r  \der z = (z_2 - z_1)\cdot 2\pi r
		.\end{align*}
\end{enumerate}

\begin{figure}[h]
	\centering
	\def\svgwidth{\columnwidth}
\begingroup%
  \makeatletter%
  \providecommand\color[2][]{%
    \errmessage{(Inkscape) Color is used for the text in Inkscape, but the package 'color.sty' is not loaded}%
    \renewcommand\color[2][]{}%
  }%
  \providecommand\transparent[1]{%
    \errmessage{(Inkscape) Transparency is used (non-zero) for the text in Inkscape, but the package 'transparent.sty' is not loaded}%
    \renewcommand\transparent[1]{}%
  }%
  \providecommand\rotatebox[2]{#2}%
  \newcommand*\fsize{\dimexpr\f@size pt\relax}%
  \newcommand*\lineheight[1]{\fontsize{\fsize}{#1\fsize}\selectfont}%
  \ifx\svgwidth\undefined%
    \setlength{\unitlength}{712.71753506bp}%
    \ifx\svgscale\undefined%
      \relax%
    \else%
      \setlength{\unitlength}{\unitlength * \real{\svgscale}}%
    \fi%
  \else%
    \setlength{\unitlength}{\svgwidth}%
  \fi%
  \global\let\svgwidth\undefined%
  \global\let\svgscale\undefined%
  \makeatother%
  \begin{picture}(1,0.40581736)%
    \lineheight{1}%
    \setlength\tabcolsep{0pt}%
    \put(0,0){\includegraphics[width=\unitlength,page=1]{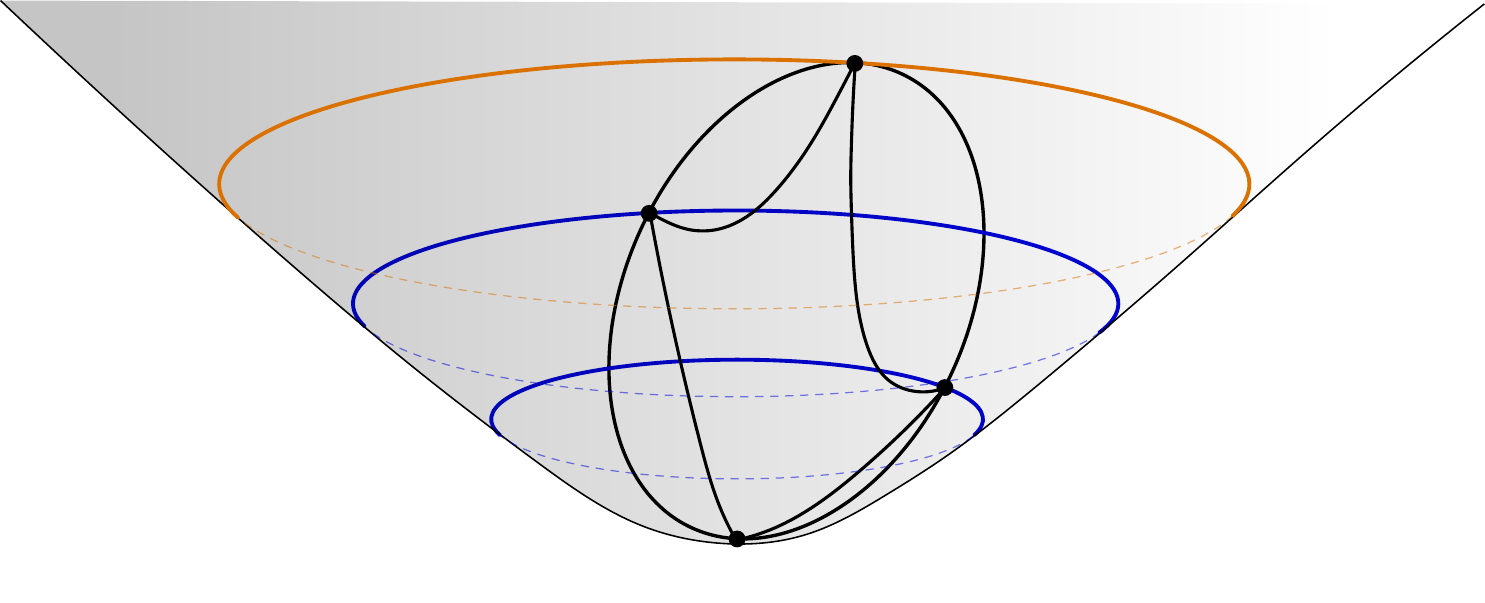}}%
    \put(0.6058089,0.37905507){\makebox(0,0)[lt]{\lineheight{0}\smash{\begin{tabular}[t]{l}$D$\end{tabular}}}}%
    \put(0.52493127,0.0016497){\makebox(0,0)[lt]{\lineheight{0}\smash{\begin{tabular}[t]{l}$A$\end{tabular}}}}%
    \put(0.65302031,0.14089834){\makebox(0,0)[lt]{\lineheight{0}\smash{\begin{tabular}[t]{l}$B$\end{tabular}}}}%
    \put(0.40442217,0.27757434){\makebox(0,0)[lt]{\lineheight{0}\smash{\begin{tabular}[t]{l}$C$\end{tabular}}}}%
    \put(0.66408066,0.07837504){\makebox(0,0)[lt]{\lineheight{0}\smash{\begin{tabular}[t]{l}$O_{AB}$\end{tabular}}}}%
    \put(0.75058608,0.1511838){\makebox(0,0)[lt]{\lineheight{0}\smash{\begin{tabular}[t]{l}$O_{AC}$\end{tabular}}}}%
    \put(0.83709151,0.22687603){\makebox(0,0)[lt]{\lineheight{0}\smash{\begin{tabular}[t]{l}$O_{AD}$\end{tabular}}}}%
  \end{picture}%
\endgroup%

	\caption{\Cref{thm:diag_pyth} on a hyperboloid}
	\label{fig:proof_hyperbolic}
\end{figure}
We will briefly repeat the proof of \Cref{thm:diag_pyth} on the hyperboloid, also see \cref{fig:proof_hyperbolic}.
By (a) we may assume that $A$ is the point $(0,0,r)$. As the points $A,B,D,C$ lie on a circle in  $\mathbb{H}^2$, (b) tells us that they lie in a Euclidean plane $P$. 
As the arcs $\arc {AC}, \arc {BD}$ are congruent, as well as $\arc {AB},\arc {CD}$, we conclude that $|AC| = |BD|, |AB| = |CD|$ from (a). 
As these points lie in a plane that is Euclidean we conclude that $ABDC$ is a parallelogram in $P$.
Thus
\[
\overrightarrow{BA}_z + \overrightarrow{CA}_z = \overrightarrow{DA}_z
,\]
and we deduce the theorem from (c). 

\begin{remark}
	We use the sphere in Euclidean space for the spherical case, and the hyperboloid in Minkowski space for the hyperbolic case. 
	In fact, we can use the same proof to also obtain the Pythagorean theorem in the flat case by working on the paraboloid $x^2 + y^2 = z$, which we identify as the Euclidean plane by pulling back all geometric information from the projection to the $xy$-plane. 
	However, this argument is cyclic as it depends on knowing the Pythagorean theorem in $\E^3$.
\end{remark}

\section{Axioms and embeddings} \label{sec:axioms_and_embeddings}
Let me preface this section by stating that I know very little about axiomatic geometry and would welcome feedback very much. 

As remarked by Pambuccian \cite{pambuccianMariaTeresaCalapsos2010}, a true axiomatic proof of \Cref{thm:diag_pyth} using Hilbert's axioms is impossible, because it is impossible to treat the area of circles. 
Of course, if we are in any setting where we form a sort of measure theory we would be close to Riemannian geometry in which the proof can obviously be worked out, as $\sphere$ has a unique embedding in $\E^3$ up to isometries of $\E^3$. 
So the question is if there is any weaker setting in which the result can be phrased and proved. 

For polygons, we can treat area essentially by defining scissor congruence or equivalently the defect of the sum of angles. But this fails for circles.  
If we have already embedded $\sphere$ into $\E^3$, we can take the bread-crust theorem as the definition of the area of disks in $\sphere \subset \E^3$, similarly to how we can define the area of polygons as defect of the angles between the planes in $\E^3$ corresponding to the arcs in $\sphere$.
In this way, it is easy to see that \Cref{thm:diag_pyth} and the proof of \cref{sec:the_actual_proof} can be done entirely synthetically in  $\E^3$.
I would like to think that if one has a good axiomatization of $\mink$, the hyperbolic case can also be treated that way, but I do not know of any such axiomatization. 

I wonder if this is not already close to the best we can do. 
The bread-crust theorem tells us that a theory of areas of disks is close to an embedding in $\E^3$. 
Indeed, fix a center point $a \in \sphere$, then for any $b \in \sphere$ the area of the disk $O_{AB}$ acts as a coordinate of $B$ in $\E^3$.  
So any theorem of 3-dimensional Euclidean geometry could be translated to a theorem about coordinates and thus about areas of disks on $\sphere$.
Therefore if we have a set of axioms that builds on Hilbert's axioms to describe the areas of disks, this would already be very close to considering the sphere in $\E^3$, similarly for the hyperboloid in $\mink$.

\section{Acknowledgements} \label{sec:acknowledgements}
At a conference in CIRM in Marseille, Patrick Popescu-Pampu introduced me to the absolute Pythagorean theorem and Paolo Maraner's paper on it \cite{maranerSphericalPythagoreanTheorem2010}. 
He mentioned that so far only an analytic proof is known, and a synthetic proof would be of interest. 

Special gratitude goes to Patrick Popescu-Pampu for introducing me into the problem, encouraging me to write this paper, suggesting that the proof might carry over to the hyperboloid model, digging up old papers whose content is difficult to find nowadays and providing feedback on the paper.
I would also like to thank CIRM for creating such a wonderful place for mathematicians to find each other and ideas to grow.
Finally, I want to thank Joeri Van der Veken and Johannes Nicaise for discussing the proof with me and providing feedback.

\section*{Declarations} \label{sec:declarations}
\textbf{Funding:} The author was supported by the long term structural funding (Methusalem grant) by the Flemish Government.

\printbibliography

@article{reynoldsHyperbolicGeometryHyperboloid1993,
    AUTHOR = {Reynolds, William F.},
     TITLE = {Hyperbolic geometry on a hyperboloid},
   JOURNAL = {Amer. Math. Monthly},
  FJOURNAL = {American Mathematical Monthly},
    VOLUME = {100},
      YEAR = {1993},
    NUMBER = {5},
     PAGES = {442--455},
      ISSN = {0002-9890,1930-0972},
   MRCLASS = {51M10},
  MRNUMBER = {1215530},
MRREVIEWER = {Hanna\ Sandler},
       DOI = {10.2307/2324297},
       %URL = {https://doi.org/10.2307/2324297},
}

@article{pambuccianMariaTeresaCalapsos2010,
  title = {Maria {{Teresa Calapso}}’s {{Hyperbolic Pythagorean Theorem}}},
  author = {Pambuccian, Victor},
  date = {2010-12-01},
  journaltitle = {The Mathematical Intelligencer},
  shortjournal = {Math Intelligencer},
  volume = {32},
  number = {4},
  pages = {2--2},
  issn = {1866-7414},
  doi = {10.1007/s00283-010-9169-0},
  langid = {english}
}

@article{maranerSphericalPythagoreanTheorem2010,
    AUTHOR = {Maraner, Paolo},
     TITLE = {A spherical {P}ythagorean theorem},
   JOURNAL = {Math. Intelligencer},
  FJOURNAL = {The Mathematical Intelligencer},
    VOLUME = {32},
      YEAR = {2010},
    NUMBER = {3},
     PAGES = {46--50},
      ISSN = {0343-6993,1866-7414},
   MRCLASS = {51M04},
  MRNUMBER = {2721310},
       DOI = {10.1007/s00283-010-9152-9},
       %URL = {https://doi.org/10.1007/s00283-010-9152-9},
}

@video{3blue1brownWhySpheresSurface2018,
  entrysubtype = {video},
  title = {But Why Is a Sphere's Surface Area Four Times Its Shadow?},
  editor = {{3Blue1Brown}},
  editortype = {director},
  date = {2018-12-02},
  url = {https://www.youtube.com/watch?v=GNcFjFmqEc8},
  urldate = {2025-02-20}
}

@article {calapso,
    AUTHOR = {Familiari-Calapso, Maria Teresa},
     TITLE = {Sur une classe de triangles et sur le th\'eor\`eme de
              {P}ythagore en g\'eom\'etrie hyperbolique},
   JOURNAL = {C. R. Acad. Sci. Paris S\'er. A-B},
  FJOURNAL = {Comptes Rendus Hebdomadaires des S\'eances de l'Acad\'emie des
              Sciences. S\'eries A et B},
    VOLUME = {268},
      YEAR = {1969},
     PAGES = {A603--A604},
      ISSN = {0151-0509},
   MRCLASS = {50.40},
  MRNUMBER = {247565},
MRREVIEWER = {F.\ Bachmann},
}

@book{thurstonThreedimensionalGeometryTopology1997,
    AUTHOR = {Thurston, William P.},
     TITLE = {Three-dimensional geometry and topology. {V}ol. 1},
    SERIES = {Princeton Mathematical Series},
    VOLUME = {35},
    EDITOR = {Levy, Silvio},
 PUBLISHER = {Princeton University Press, Princeton, NJ},
      YEAR = {1997},
     PAGES = {x+311},
      ISBN = {0-691-08304-5},
   MRCLASS = {57M50 (53A35 57M25 57M60 57N10)},
  MRNUMBER = {1435975},
MRREVIEWER = {Athanase\ Papadopoulos},
}
\end{document}